\documentclass{article}
\usepackage[utf8]{inputenc}
\usepackage[usenames]{color}
\usepackage{amsmath}
\usepackage{epsfig,latexsym}
\usepackage{epsf}
\usepackage{amsmath,amsfonts,amstext,amssymb,amsbsy,amsopn,amsthm,eucal}
\usepackage{latexsym}
\usepackage{epsfig,latexsym}
\usepackage[usenames]{color}
\usepackage{dirtytalk}

\usepackage{indentfirst}

\title{Ouroboros Functionals, Families of Ouroboros Functions, and Their Relationship to Partial Differential Equations and Probability Theory}
\author{Nathan Thomas Provost$\footnote{Student of Applied Mathematics and Statistics at Brown University. University Email: nathan\_provost@brown.edu}$}
\date{}

\begin{document}

\maketitle

\begin{center}
\textbf{Abstract}
\end{center}

\small Previously, we have introduced a very small number of examples of what we call Ouroboros functions. Using our already established theory of Ouroboros spaces and their functions, we will provide a set of families of Ouroboros functions that bolster our overall understanding of the Ouroboros spaces. From here, we extend the theory of Ouroboros functions by introducing Ouroboros functionals and Ouroboros functional spaces. Furthermore, we re-frame the expected value of a random variable as an Ouroboros functional, which proves to be more intuitive in view of probabilistic measure theory. We then show that these Ouroboros functions have additional applications, as they are general solutions to certain elementary linear first order partial differential equations (PDEs). We conclude by elaborating upon this connection and discussing future endeavors, which will be centered on answering a given hypothesis.

\normalsize

\small

\section*{Introduction}

As in previous papers (\cite{1} \cite{2}), we are concerned with the Ouroboros spaces, whose definition is: $\textbf{\textit{O}}(A^n)=\left\{f:A^n\rightarrow B \mid f(f(\textbf{x}),...,f(\textbf{x}))=f(\textbf{x}), \ \forall \textbf{x}\in A^n, \ \forall B\subseteq A\right\}.$ We have appropriately referred to the functions contained therein as Ouroboros functions. However, we have not presented many cases of these functions. On the contrary, we have previously only given functions of the forms:
\[f(x_1,...,x_n)=\frac{1}{n}\sum_{i=1}^{n}x_i \ \ \rm{or} \ \ \textit{$f(x_1$},...,\textit{$x_n)$}=\textit{c}\in\mathbb{R}\]
This only provides us with two Ouroboros functions per domain of dimension \textit{n}, which is quite limited in retrospect since, for example, we can come up with infinitely many quadratic functions or cubic functions. Yet recently, we have constructed several infinite families of nontrivial Ouroboros functions that branch off of this arithmetic average function, which turns out to be a special case of a more general form, which also solves certain linear partial differential equations (PDEs). Furthermore, we can extend the idea of the Ouroboros functions to functionals, which arises from the acknowledgement of these general forms.

\newpage

\normalsize

\section*{Constructing the Ouroboros Function Families}

These function families arose from our examination of two specific linear first order PDEs, which we will discuss later. However, for explanatory purposes, we will begin by referencing a familiar example of an Ouroboros function. We have previously shown that $f(x,y)=\frac{1}{2}(x+y)$ is an Ouroboros function; more specifically, we showed that $f\in\textbf{\textit{O}}(\mathbb{R}^2)$ (see \cite{1}). We can write $f$ in the form $f(x,y)=\frac{1}{2}x+\frac{1}{2}y=c_1x+c_2y$, where $c_1=\frac{1}{2}=c_2$. Two conditions hold in this case for $c_1$ and $c_2$: $c_1+c_2=1$ and $c_1=c_2$. We wish to preserve the first condition, but we will now try to relax the second.\\

Suppose now that $c_1=\frac{1}{4}$ and $c_2=\frac{3}{4}$, such that $c_1+c_2=1$, but $c_1\neq c_2$. Our function is now of the form $f(x,y)=f=\frac{1}{4}x+\frac{3}{4}y$. Naturally, $f(f,f)=\frac{1}{4}\left(\frac{1}{4}x+\frac{3}{4}y\right)+\frac{3}{4}\left(\frac{1}{4}x+\frac{3}{4}y\right)=\frac{1}{4}\left(\frac{1}{4}x+\frac{3}{4}y+3\left(\frac{1}{4}x+\frac{3}{4}y\right)\right)=\frac{1}{4}\left(\frac{1}{4}x+\frac{3}{4}y+\frac{3}{4}x+\frac{9}{4}y\right)=\frac{1}{4}\left(x+3y\right)=\frac{1}{4}x+\frac{3}{4}y=f$. Also, if we assume $x,y\in\mathbb{R}$, as we did with our prior function, $f:\mathbb{R}^2\rightarrow\mathbb{R}$, which cumulatively means that $f\in\textbf{\textit{O}}(\mathbb{R}^2)$. From this observational case, we move on and prove the following general theorem.\\

\textbf{Theorem}: Let $c_i\in\mathbb{R}, \ \forall i\in\{1,...,n\}\subset\mathbb{N}$ and assume $x_1,...,x_n\in\mathbb{R}$. Then:
\[\sum_{i=1}^{n}c_i=1 \longrightarrow f(x_1,...,x_n)=\sum_{i=1}^{n}c_ix_i\in\textbf{\textit{O}}(\mathbb{R}^n)\]

\noindent\textbf{Proof}: Assume that the appropriate conditions are met, most importantly that the coefficients given by $c_i$ are real numbers that add up to one. Letting $f(x_1,...,x_n)=f$, we see that:
\[f(f,...,f)=c_1\left(\sum_{i=1}^{n}c_ix_i\right)+...+c_n\left(\sum_{i=1}^{n}c_ix_i\right)=(c_1+...+c_n)\sum_{i=1}^{n}c_ix_i=\]
\[\sum_{i=1}^{n}c_i\sum_{i=1}^{n}c_ix_i=\sum_{i=1}^{n}c_ix_i=f, \ \rm{since} \ \sum_{\textit{i}=1}^{\textit{n}}\textit{c}_\textit{i}=1\]
Similarly, we know that $x_1,...,x_n$ are all real-valued variables, which means that $f:\mathbb{R}^n\rightarrow\mathbb{R}$ and $f(f,...,f)=f$. Therefore, $f\in\textbf{\textit{O}}(\mathbb{R}^n)$. $\qedsymbol$\\

We see that for each domain of dimension $n$, an infinite family of Ouroboros functions is formed. The case mentioned in our previous discussion (\cite{1}) was a special example, where $c_1=...=\frac{1}{n}=...=c_n$. This important special case is simply the arithmetic average of our set of real variables. While it has some nice properties (as mentioned in \cite{1} and discussed later), these general Ouroboros functions can help give us an idea of what an Ouroboros function looks like. Even so, these families are limited, since they are all linear, with each variable having an exponent that does not exceed 1.\\

\section*{Expected Value as an Ouroboros Functional}

From here, we notice that a more intricate relationship between probability theory and the Ouroboros spaces exists. Suppose that we impose another condition on the constants, namely that:
\[\sum_{i=1}^{n}c_i=1 \ \rm{and} \ \textit{c}_\textit{i}\in[0,1], \ \forall\textit{i}\in\{1,...,\textit{n}\}\subset\mathbb{N}\]
If we define $c_i=\mathbb{P}(X=x_i)$ for a discrete, simple random variable on a probability space $\mathbb{S}=(\Omega,\mathbb{F},\mathbb{P})$ and treat $x_1,...,x_n$ as the values taken on by $X$ instead of variables, then we see that our function is now given by:
\[f(x_1,...,x_n)=\sum_{i=1}^{n}x_i\mathbb{P}(X=x_i)=\mathbb{E}[X]\]
This definition corresponds to the definition of expected value as given in traditional sources (such as \cite{3} \cite{4} \cite{5}). It is important to note again, however, that the $x_i$ values are no longer variables in this case, but are instead the singular values that $X=X(\omega)$ takes on, since $X:\Omega\rightarrow\{x_1,...,x_n\}\subset\mathbb{R}$, $\forall\omega\in\Omega$ in this case. Accordingly, we can think of the expected value of a random variable $(\mathbb{E}[X]=\mathbb{E})$ for the space of all real-valued random variables given by $\mathbb{X}_{\mathbb{R}}$ as a functional given by $\mathbb{E}:\mathbb{X}_{\mathbb{R}}\rightarrow\mathbb{R}$. We therefore offer the following definition, which is a functional analog to the Ouroboros spaces.\\

\textbf{Definition}: The \textbf{Ouroboros Functional Space} for the function space given by $\mathcal{F}_\mathbb{D}$ (whose functions either have domain $\mathbb{D}$ or have a domain contained in $\mathbb{D}$) is the set of functionals given by:
\[\tilde{\textbf{\textit{O}}}(\mathcal{F}_\mathbb{D})=\{\mathbb{O}:\mathcal{F}_\mathbb{D}\rightarrow B \ | \ \mathbb{O}(\mathbb{O}(f))=\mathbb{O}(f), \ \forall f\in\mathcal{F}_\mathbb{D}, \ \forall B\subseteq\mathbb{D}\}\]
An element of this space is called an \textbf{Ouroboros functional}. Furthermore, for a functional that takes several functions from the same function space as inputs, we can generalize this definition to:
\[\tilde{\textbf{\textit{O}}}(\mathcal{F}^n_\mathbb{D})=\{\mathbb{O}:\mathcal{F}^n_\mathbb{D}\rightarrow B \ | \ \mathbb{O}(\mathbb{O}(\textbf{f}),...,\mathbb{O}(\textbf{f}))=\mathbb{O}(\textbf{f},...,\textbf{f}), \ \forall \textbf{f}\in\mathcal{F}^n_\mathbb{D}, \ \forall B\subseteq\mathbb{D}\}\]
Finally, we can consider a mixed Cartesian product of function spaces given by $\mathcal{P}=\mathcal{F}_{\mathbb{D}_1}\times...\times\mathcal{F}_{\mathbb{D}_n}$, for which the Ouroboros functional space would be:
\[\tilde{\textbf{\textit{O}}}(\mathcal{P})=\{\mathbb{O}:\mathcal{P}\rightarrow B \ | \ \mathbb{O}(\mathbb{O}(\textbf{f}),...,\mathbb{O}(\textbf{f}))=\mathbb{O}(\textbf{f},...,\textbf{f}), \ \forall \textbf{f}\in\mathcal{P}, \ \forall B\subseteq\mathbb{D}_i, \forall i=1,...,n\}\]
The reason we stated that the functions of $\mathcal{F}_{\mathbb{D}}$ have domains that are contained in or equal to $\mathbb{D}$ stems from the fact that the commonly used number systems are subsets of one another ($\mathbb{N}\subset\mathbb{Z}\subset\mathbb{Q}\subset\mathbb{R}$ etc.) Therefore, a function whose domain is the natural numbers is also technically a function whose domain is the real numbers. We have already discussed a loose interpretation of $\mathbb{E}[X]$ as an Ouroboros functional, but now we should formally justify this claim.

\newpage

\textbf{Proposition 1}: Let $\mathbb{X}_{\mathbb{R}}$ be the space of real-valued random variables for the probability spaces given by $\mathbb{S}_X=(\Omega_X,\mathbb{F}_X,\mathbb{P}_X)$, $\forall X\in\mathbb{X}_\mathbb{R}$. Assume that for all $X\in\mathbb{X}_{\mathbb{R}}$:
\[\mathbb{E}[X]=c_X\in\mathbb{R} \ \rm{and} \ \int_{\Omega_X}|\textit{X}|\textit{d}\mathbb{P}_{\textit{X}}<\infty \]
Then $\mathbb{E}:\mathbb{X}_{\mathbb{R}}\rightarrow\mathbb{R} \ni \mathbb{E}\in\tilde{\textbf{\textit{O}}}(\mathbb{X}_\mathbb{R})$.\\

\noindent\textbf{Proof}: We accept the aforementioned assumptions and observe, as noted, that $\mathbb{E}:\mathbb{X}_\mathbb{R}\rightarrow\mathbb{R}$ such that $\mathbb{R}\subseteq\mathbb{R}$. As noted in numerous sources (like \cite{4} and \cite{5}), the definition of expected value for some random variable $X\in\mathbb{X}_\mathbb{R}$ for the probability space $\mathbb{S}_X=(\Omega_X,\mathbb{F}_X,\mathbb{P}_X)$ is given by:
\[\mathbb{E}[X]=\int_{\Omega_X}Xd\mathbb{P}_X\]
We can define a simple random variable $E$ such that $\mathbb{E}[X]=E=c_X=c_X\textbf{1}_{\{\omega_X\in\Omega_X\}}$, where:
\[\textbf{1}_{\{\omega_X\in\Omega_X\}}=\begin{cases}
1 & \omega_X\in\Omega_X \\
0 & \omega_X\notin\Omega_X
\end{cases}\]
Using the definition of the Lebesgue integral for simple random variables, we see that:
\[\mathbb{E}[\mathbb{E}[X]]=\mathbb{E}[E]=\int_{\Omega_X}c_X\textbf{1}_{\{\omega_X\in\Omega_X\}}d\mathbb{P}_X=c_X\mathbb{P}_X(\{\omega_X\in\Omega_X\})\]
By the definition of a probability measure for a probability space (again, as defined in sources like \cite{4} and \cite{5}), we know that $\mathbb{P}_X({\{\omega_X\in\Omega_X\}})=1$, so $c_X\mathbb{P}_X({\{\omega_X\in\Omega_X\}})=c_X=\mathbb{E}[X]$. Therefore, $\mathbb{E}[\mathbb{E}[X]]=\mathbb{E}[X], \forall X\in\mathbb{X}_\mathbb{R}$, which cumulatively means that $\mathbb{E}\in\tilde{\textbf{\textit{O}}}(\mathbb{X}_\mathbb{R})$. $\qedsymbol$\\

\noindent \textbf{Remark}: Previously, we have shown that (under a set of \textit{specific conditions}) $\mathbb{E}[X]\in\textbf{\textit{O}}(\mathbb{R}^{\infty})$ almost surely, making $\mathbb{E}$ an Ouroboros function for an infinite domain $(\mathbb{R}^\infty)$. This, however, is a counter-intuitive characterization of the expected value functional, though it is technically correct under the proper circumstances (see \cite{1}). Nonetheless, it is more appropriate to adopt the convention of treating the expected value of a random variable as an Ouroboros functional, as it applies more generally for all real-valued variables, whereas the previous characterization only holds for random variables whose domain is $\mathbb{R}$, and not a subset of $\mathbb{R}$.\\

The definition of Ouroboros functionals is simply the next logical step in the process of extending the theory of Ouroboros spaces. We are now able to divide functionals (specifically higher-order functions) into analogous Ouroboros functional spaces. This provides a more sound basis for the classification of a random variable's expected value, since being an Ouroboros functional allows for greater generalization and intuition when considering general random variables.

\newpage

\section*{Ouroboros Functions as Solutions to Linear  First Order Partial Differential Equations}

Initially, our reconsideration of Ouroboros spaces began in the context of investigating alternative, aesthetic solutions to linear first order partial differential equations. Our exploration, while certainly not overtly elaborate in nature, revealed some particularly interesting relationships between two linear first order PDEs and certain families of Ouroboros functions. The two PDEs in question are as follows:
\[\sum_{k=1}^{\beta}\frac{\partial u}{\partial x_k}=\beta\frac{\partial u}{\partial x_n}, \ \forall\beta\in\mathbb{N}\ni1\leq\beta\leq n \ \ \ \rm{(I)}\]
\[\sum_{k=1}^n(-1)^k\frac{\partial u}{\partial x_k}=0 \ \ \ \rm{(II)}\]
We first examine these equations individually, and then proceed to consider them as a system. Both of these equations are transport PDEs, which have well-documented solutions, usually derived through the method of characteristics (as described in texts like \cite{6} and \cite{7}). However, we aimed to take a looser, but still mathematically correct approach to these equations with the hope of connecting them to Ouroboros functions and showcasing some of their more aesthetic solutions.\\

We begin by considering equation (I). A linear function with constant coefficients appears to be a fairly simple candidate, and so we will construct a solution of this form. First, let us consider the case where $\beta=n$. If the only term of the function $u(x_1,...,x_n)$ containing $x_n$ is of the form $c_nx_n$ for some $c_n\in\mathbb{R}$, let us choose $c_n=\frac{1}{n}$. Similarly, suppose $u$ is a linear combination of all its variables. More specifically, let $u$ be of the form:
\[u(x_1,...,x_n)=\sum_{k=1}^{n}c_kx_k\]
where $c_1,...,c_n\in\mathbb{R}$ and $c_n=\frac{1}{n}$. Equation (I) is now reduced to the following:
\[\sum_{k=1}^{n}\frac{\partial u}{\partial x_k}=1 \rightarrow \sum_{k=1}^{n}c_k=1\]
Evidently, it holds that:
\[\mathbb{U}=\left\{u:\mathbb{R}^n\rightarrow\mathbb{R} \ | \ u(x_1,...,x_n)=\sum_{k=1}^{n}c_kx_k, \ \forall c_k\in\mathbb{R} \ni c_n=\frac{1}{n}, \ \sum_{k=1}^{n}c_k=1 \right\}\]
consists solely of functions that solve equation (I) for $\beta=n$. Similarly, following from the theorem we previously proved, it holds that $\mathbb{U}\subset\textbf{\textit{O}}(\mathbb{R}^n)$, since all functions with the same form as $u$ whose constant coefficients add up to 1 are Ouroboros functions.

\newpage

For the case where $\beta\in\{1,...,n-1\}$, we will now verify the following proposition regarding the solution to equation (I).\\

\textbf{Proposition 2}: Any function of the form:
\[u(x_1,...,x_n)=\left(\left[\frac{1}{\beta}\sum_{k=1}^{\beta}c_k\right]x_n+\sum_{k=1}^{n-1}c_kx_k\right) \ \ \ \begin{cases}
\forall \beta\in\{1,...,n-1\}\subset\mathbb{N}\\
\forall c_k\in\mathbb{R}, \ \forall k\in\{1,...,n-1\}\subset\mathbb{N}
\end{cases}\]
solves equation (I) under the previously stated conditions next to $u$.\\

\noindent\textbf{Proof}: Assume $u$ is given as above and all of the other parameters fall in line with the specifications of the proposition. Naturally, we see that:
\[\frac{\partial u}{\partial x_n}=\frac{1}{\beta}\sum_{k=1}^{\beta}c_k\ \therefore \ \beta\frac{\partial u}{\partial x_n}=\sum_{k=1}^{\beta}c_k\]
Conclusively, we further note that:
\[\frac{\partial u}{\partial x_k}=c_k, \ \forall k< n \ \therefore \ \sum_{k=1}^{\beta}\frac{\partial u}{\partial x_k}=\sum_{k=1}^{\beta}c_k=\beta\frac{\partial u}{\partial x_n} \ \ \ \qedsymbol\]
We further note that the coefficient of $x_n$ is simply the arithmetic average of the coefficients $c_1,...,c_{\beta}$. Combining the previous result, we may wish to rewrite $u$, as a general solution to equation (I), in the form:
\[u(x_1,...,x_n)=\mu_{\beta}x_n+\sum_{k=1}^{n-1}c_kx_k \ \ni \ \mu_{\beta}=\frac{1}{\beta}\sum_{k=1}^{\beta}c_k \ \ \rm{if} \ \ \begin{cases}
\beta\in\{1,...,n-1\}\subset\mathbb{N}\\
c_k\in\mathbb{R}, \ \forall k\in\{1,...,n\}\subset\mathbb{N}
\end{cases}\]
\[u(x_1,...,x_n)=\sum_{k=1}^{n}c_kx_k \ \rm{if} \ \forall \textit{c}_\textit{k}\in\mathbb{R} \ni \textit{c}_\textit{n}=\frac{1}{\textit{n}}, \ \sum_{\textit{k}=1}^{\textit{n}}\textit{c}_\textit{k}=1 \ \rm{and} \ \beta=\textit{n}\]
While it is possible to obtain even more general solutions to equation (I) through the method of characteristics, these solutions show that there exists an aesthetic connection between Ouroboros spaces and the solutions to equation (I). Furthermore, we can see similar connections by exploring equation (II). Again, suppose that our solution takes the form of a linear combination of each variable multiplied by a constant coefficient, which must be a real number.\\

\textbf{Proposition 3}: If $n\in\mathbb{N}$ is even, then:
\[u(x_1,...,x_n)=\sum_{k=1}^{n}c_kx_k \ \rm{is \ a \ solution \ to \ equation \ (II) \ if} \ \sum_{\textit{k}=1,3,5}^{\textit{n}-1}\textit{c}_\textit{k}=\sum_{\textit{k}=2,4,6}^{\textit{n}}\textit{c}_\textit{k}\]

\noindent\textbf{Proof}: Suppose $u$ is a function as given above and $n$ is an even natural number. By basic algebraic properties, the following holds for equation (II):
\[\sum_{k=1}^n(-1)^k\frac{\partial u}{\partial x_k}=0 \iff\sum_{k=2,4,6}^{n}\frac{\partial u}{\partial x_k}-\sum_{k=1,3,5}^{n-1}\frac{\partial u}{\partial x_k}=0 \ \therefore \ \sum_{k=2,4,6}^{n}\frac{\partial u}{\partial x_k}=\sum_{k=1,3,5}^{n-1}\frac{\partial u}{\partial x_k}\]
\newpage

Since $u$ is given by a linear combination of its variables with constant coefficients, it holds that:
\[\sum_{k=2,4,6}^{n}\frac{\partial u}{\partial x_k}=\sum_{k=2,4,6}^{n}c_k \ \rm{and} \ \sum_{\textit{k}=1,3,5}^{\textit{n}-1}\frac{\partial \textit{u}}{\partial \textit{x}_\textit{k}}=\sum_{\textit{k}=1,3,5}^{\textit{n}-1}\textit{c}_\textit{k}\]
Therefore, it holds that $u(x_1,...,x_n)$ as given is a solution to equation (II) if:
\[\sum_{k=2,4,6}^{n}c_k=\sum_{k=1,3,5}^{n-1}c_k \ \ \ \qedsymbol\]
Naturally, if all of the constants are equal to one another, the previously solution still holds if $n$ is even. This fact becomes important when we attempt to solve equations (I) and (II) simultaneously (hereinafter referred to as \say{the system}). For the purpose of our study, we will consider the system under the condition $\beta=n$ in equation (I), where $n$ is an even natural number. Additionally, we are interested in imposing an initial condition (given below) on the system. Subsequently, we now prove the following proposition.\\

\textbf{Proposition 4}: The function given by:
\[u(x_1,...,x_n)=\frac{1}{n}\sum_{k=1}^{n}x_k, \ \forall x_k\in\mathbb{R}, \ \forall k\in\{1,...,n\}\]
solves the system of differential equations given by:
\[\sum_{k=1}^{n}\frac{\partial u}{\partial x_k}=n\frac{\partial u}{\partial x_n}; \ \forall n\in\mathbb{N}\ni \rm{\textit{n} \ is \ even}\]
\[\sum_{k=1}^n(-1)^k\frac{\partial u}{\partial x_k}=0; \  u(0,...,0)=0\]\\

\noindent\textbf{Proof}: Let our solution $u(x_1,...,x_n)=u$ for an even number $n$ be given by:
\[u(x_1,...,x_n)=\frac{1}{n}\sum_{k=1}^{n}x_k\]
From this definition, we make the following observations:
\[\frac{\partial u}{\partial x_n}=\frac{1}{n} \ni n\frac{\partial u}{\partial x_n}=1, \ \rm{and \ furthermore,} \ \frac{\partial \textit{u}}{\partial \textit{x}_\textit{k}}=\frac{1}{\textit{n}}, \forall \textit{k}\in\{1,...,\textit{n} \ \therefore\}\]
\[\sum_{k=1}^{n}\frac{\partial \textit{u}}{\partial \textit{x}_\textit{k}}=\sum_{k=1}^{n}\frac{1}{n}=1=n\frac{\partial u}{\partial x_n}\]

\newpage

Now that we know that $u$ solves equation (I), we see that $u$ is of the form:
\[u(x_1,...,x_n)=\sum_{k=1}^{n}c_kx_k \ \rm{and} \ \sum_{\textit{k}=1,3,5}^{\textit{n}-1}\textit{c}_\textit{k}=\sum_{\textit{k}=1,3,5}^{\textit{n}-1}\frac{1}{\textit{n}}=\frac{\textit{n}}{2}=\sum_{\textit{k}=2,4,6}^{\textit{n}}\frac{1}{\textit{n}}=\sum_{\textit{k}=2,4,6}^{\textit{n}}\textit{c}_\textit{k}\]
Therefore, $u$ solves equation (II) by Proposition 3. Finally, it is evident that:
\[u(0,...,0)=\frac{1}{n}\sum_{k=1}^{n}0=0\]
Hence, we conclude that $u$ solves the given system of PDEs. $\qedsymbol$

\section*{Conclusion}

The aesthetic intrigue that arises from our proof of Proposition 4 originates from the fact that $u$ is an Ouroboros function belonging to $\textbf{\textit{O}}(\mathbb{R}^n)$. While the PDEs we have previously mentioned are generally well studied (as noted in \cite{6} and \cite{7}), our investigation into the ties between the two equations and Ouroboros functions has shown potential for further research. On a more general note, the arithmetic average as a concept seems to be deeply tied to these two transport PDEs. The coefficient of $x_n$ in our set of solutions to equation (I) for $\beta\in\{1,...,n-1\}$ is the arithmetic average of the coefficients $c_1,...,c_\beta$, and the general arithmetic average function is not only an Ouroboros function, but also a solution to both equations (I) and (II) (including the initial condition). The probabilistic idea of an average (\textit{i}.\textit{e}. the expected value of a random variable) is also intricately related to the general theory of Ouroboros functions, just as we showed in Proposition 1. We demonstrated that it is more intuitive to think of the expected value of a random variable as an Ouroboros \textit{functional}, rather than an Ouroboros function, especially in the broader context of probabilistic measure theory. Moreover, the generalization of the Ouroboros spaces in the form of the Ouroboros functional spaces will undoubtedly allow for a more thorough investigation of these concepts in the future. Finally, Proposition 4 has laid the foundation for an area of future investigations. Our future research efforts will be centered on proving or disproving the following hypothesis.\\

\textbf{Hypothesis}: If
\[\sum_{k=1}^{n}\frac{\partial u}{\partial x_k}=n\frac{\partial u}{\partial x_n}; \ \forall n\in\mathbb{N}\ni \rm{\textit{n} \ is \ even}\]
\[\sum_{k=1}^n(-1)^k\frac{\partial u}{\partial x_k}=0; \ \ u\in\textbf{\textit{O}}(\mathbb{R}^n),\]
then
\[u(x_1,...,x_n)=\frac{1}{n}\sum_{k=1}^{n}x_k, \ \forall x_k\in\mathbb{R}, \ \forall k\in\{1,..,n\}\]
is the \textit{unique}, \textit{nontrivial} solution to this overdetermined system of PDEs.

\end{document}